\documentclass[a4paper]{amsart}  

\usepackage{amssymb} \usepackage{amscd} \usepackage{times}

\usepackage[french]{babel}

\def\zbb{\mathbb{Z}}  
  
  \def\phi{\varphi}
 \def\p1{{\mathbb{P}^1_\zbb}}

\newcommand{\be} {\begin{equation}}
\newcommand{\ee} {\end{equation}}

\begin{document}

\title{Note sur les in\'egalit\'es isop\'erm\'etriques en dimension 2. }
\author{Samy Skander Bahoura}

\maketitle
\begin{abstract}
Ce texte parle de l'article de C. Bandle(1976) et A. Huber, il concerne les inegalites isoperimetriques en dimension 2 dans des cas particulier, celles de Nehari, Huber et Alexandrov. On peut trouver cela dans le livre de C. Bandle.
\end{abstract}

Soit $ D $ un domaine simplement connexe born\'e de bord $ \partial D $ analytique. En utilisant le th\'eoreme de representation conforme de Riemann, on peut supposer que $ D = B_1 $.

\bigskip

En effet, soit $ f $ une application conforme de $ D $ vers $ B_1 $, le fait que $ \partial D $ soit analytique implique qu'on peut prolonger $ f $ en application continue et injective sur $ \bar D $ tel que $ f(\partial D) = \partial B_1 $. Sans nuire a la generalite, on peut supposer que localement $ \partial D =]-1, 1[ $ et remplacer $ B_1 $ par $ B_1^+ $ une demi-boule. Dans ce cas, on a $ f(]-1, 1[)=]-1, 1[ $ et on utilise la symmetrisation de Schwarz, la continuite de $ f $ et theoreme de Morera pour prolonger $ f $ une application conforme sur un voisinage de $ D $, en considerant:

$$  g(z) = \begin{cases}
        f(z) \,\, \text{si}\,\, z\in B_1^+ \\
        \bar f({\bar z})  \,\, \text{si} \,\, z\in B_1^-
\end{cases} $$

\underbar {Etape 1:}

\bigskip

On se place sur la boule unite $ B_1 $ de $ {\mathbb R}^2 $ et $ h $ une fonction harmonique sur $ B_1 $ et qui vaut $ u $ sur le bord. $ h \in C^{\infty}(\bar B_1) $.

\bigskip

On resout l'equation: $ g $ holomorphe, $ g \in O(B_1) $,

\be |g'(z)|^2 = e^h \ee

Ceci revient a r\'esoudre par s\'eries entieres et par series de Fourier (sur le bord) une equation du type:

\be \partial \tilde g = e^{\int \partial h}, \ee

ou $ \int \partial h $ la primitive de $ \partial h $.

Il faut utiliser par exemple Parseval pour avoir la sommabilit\'e des coefficients de Fourier de $ e^{\int \partial h} $ ainsi que ses deriv\'ees successives, puis appliquer Cauchy-Schwarz pour avoir la sommabilit\'e des coefficients de $ \tilde g $ (par exemple).

(Utiliser le fait suivant, $ u= 2\log |\tilde g'(z)|-h $ est r\'eelle et remarquer que $ \partial u= \bar \partial u = 0 $).

\underbar {Etape 2:}

\bigskip

Pour l'in\'egalit\'e isop\'ermetrique de Nehari, on \'ecrit en utilisant la formule de Stokes et l'in\'egalit\'e isop\'erim\'etrique:

\be \int_{B_1} |g'(z)|^2 dx = \dfrac{1}{4} \int_{\partial B_1} g'\bar g dz-\bar g' gd\bar z=\dfrac{1}{4} \int_{g(\partial B_1)} wd\bar w-\bar wdw \ee

$$ =({\rm \,\, analogiquement,\,\, l'aire \,\,delimitee \,\, par\,\,}  g(\partial B_1) )\leq \dfrac{1}{4\pi }  {l(g(\partial B_1))}^2, $$

La preuve precedente est valable pour tout domaine avec bord regulier.

Dans tout cela, on aura utiliser la preuve par series de Fourier de l'inegalite isoperimetrique (voir le livre de B. Dacorogna sur l'inegalite isoperimetrique classique).

1- le bord d'un domaine regulier est une variete lisse de dimension 1, pour appliquer le preuve de l'inegalite isoperimetrique, il faut parametrer globalement la courbe, cette parametrisation existe (car on a un homeomorphisme vers la cercle et a partir de cet homeomorphisme  on a un diffeomorphisme lisse, par Moise-Rado, voir le Hebey,  a partir d'un homeomorphisme entre varietes lisses de dimension $ \leq 3 $ on a un diffeomorphisme lisse) 

2- Dans la preuve de l'inegalite isoperimetrique, on utilise la parmetrisation par la longeur d'arc, or ce changement de variable est $ C^1 $ seulement (car on prend la norme), par des theoremes concernant les fonctions d'une variable reelle, l'inverse est $ C^1 $, dans ce cas tout est $ C^1 $, on peut appliquer la preuve qui se base sur le theoreme de Parseval)

3- ici, soit $ z(\theta) $ une parametrisation  de la courbe definissant le bord et $ \phi (s) = \int_0^s |z'(\theta)| d \theta $ la parmetrisation par la longueur d'arc,  alors $ z'(\theta) \not = 0 $ pour tout  $ \theta $, donc, la parametrisation par la longueur d'arc $ \phi $ est lisse comme son inverse est $ C^1 $, son inverse est aussi lisse. 

On a:

$$ (\phi^{-1})'(s)=\dfrac{1}{\phi 'o \phi^{-1}(s)} $$

est lisse.

Pour le cas de l'egalite (on suppose que c'est sur la boule unite), on a le fait que $ g $ soit holomorphe $ g'\not = 0 $, $ g $ est ouverte et:

$$ \partial^* g(B_1) \subset \partial g(B_1) \subset g(\partial B_1) = C_1, $$

on en deduit  du fait que $ g(B_1) $ est ouvert, que celui-ci est un disque. Sans nuire a la generalite, on peut supposer que $ g(B_1) = B_1 $. Soit $ z_0 \in B_1 $ tel que $ g(z_0)=0 $. En composant a droite par une homographie, on peut supposer $ z_0 =0 $.

\smallskip

 On developpe $ g $ en serie entiere et on utilise le fait que $ |g(e^{i\theta})| = 1 $ (2eme  et dernier termes de $ (3) $)  pour avoir:

$$ g(z)=\sum_{n\geq 1} a_n z^n,  \,\, g(e^{i\theta})=\sum_{n\geq 1} a_n e^{i n \theta},  \,\, |z|\leq 1, \,\, \theta \in (0, 2\pi). $$

Et,

$$ \dfrac{4\pi}{4} \int_{\partial B_1} g'\bar g dz-\bar g' gd\bar z= 4\pi^2 \sum_{n\geq 1} n|a_n|^2, $$

et,

$$ 4\pi^2=2\pi \int_0^{2\pi} |g(e^{i\theta})|^2 d\theta = 4\pi^2\sum_{n\geq 1}|a_n|^2. $$

d'ou

\be \sum_{n\geq 1} n|a_n|^2 \leq \sum_{n\geq 1}|a_n|^2, \ee
d'ou, 

$$ a_n = 0, \,\, n \geq 2, \,\,\, g(z)= a_1 z, \,\, {\rm avec }\,\, |a_1|=1. $$

$ g $ est alors univalente.

\bigskip

\underbar {Etape 3:}

\bigskip

Si on pose:

$$ v= u-h, \,\, a(\mu)=\int_{ \{ v>\mu \}} |g'(z)|^2 dx, $$

Alors $ v $ verifie :

$$ \begin{cases}
 \Delta v = |g'|^2Ve^v \,\, \text{dans}\,\, B_1\\

 V\leq \lambda \\

v=0 \,\, \text{sur}\,\, \partial B_1 
\end{cases} $$

$ v $ est continue donc $ a $ est strictement decroissante, et les lignes de niveaux(grace a la formule de la coaire) sont de mesures nulles donc $ a $ est strictement decroissante et continue. De plus, en utilisant l'inegalite de Brunn-Minkowski $ a $ est localement lipschitzienne.

En effet, on utilise la formule de la coaire pour ecrire que:

$$ \int_{\{ v=\mu \}} |\nabla v| dx = \int_0^{+ \infty} H^{n-1}(\{ v= t\} \cap \{ v=\mu \}) dt = 0, $$

Ce qu'on peut ecrire:

$$  \nabla v = 0 \,\,\, \,\, \text{dx-p.p. \,\, sur }\,\,  \{ v=\mu \}, $$

On applique cela au gradient:

$$  \nabla^2 v = 0 \,\,\, \,\, \text{dx-p.p. \,\, sur }\,\,  \{ \nabla v=0 \}, $$

Donc:

$$  \nabla^2 v = 0 \,\,\, \,\, \text{dx-p.p. \,\, sur }\,\,  \{ v=\mu \}, $$

En particulier,

$$  \Delta v = 0 \,\,\, \,\, \text{dx-p.p. \,\, sur }\,\,  \{ v=\mu \}. $$

On voit alors, que si $ V >0 $ dans l'equation verifiee par $ v $, alors les lignes de niveaux sont de mesure nulles. 

\bigskip

Soient $ \mu_1 >\mu_2 $, $ D_1=\{v > \mu_1 \}, D_2=\{ v>\mu_2 \} $ et $ D_1^*, D_2^* $ les symetrises de Schwarz de $ D_1 $ et $ D_2 $, $ r_1 $ et $ r_2 $ leurs rayons respectifs, $ d=d(\partial D_1, \partial D_2) $. Grace a la formule de la moyenne, l'inegalite de Brunn-Minkowski, on obtient:

$$ \pi^{1/2} r_2= |D_2^*|^{1/2} =|D_2|^{1/2} \geq |D_1+B(0,d)|^{1/2} \geq |D_1|^{1/2}+|B(0,d)|^{1/2} $$

$$ =|D_1^*|^{1/2}+\pi^{1/2}d=\pi^{1/2}(r_1+d), $$

Comme $ v $ est continue $ \bar D_1 \subset \{ v\geq \mu_1\} $ et  $ \bar D_2 \subset \{ v \geq  \mu_2 \} $, d'ou,

$$ d(x_1, x_2)= d(\partial D_1, \partial D_2), \,\, x_1\in \partial D_1 \subset \{v=\mu_1\}, \,\, x_2 \in \partial D_2 \subset \{ v=\mu_2 \}, $$

$$ |a(\mu_1)-a(\mu_2)|=cte(||\{ v>\mu_2 \}-|\{v>\mu_1 \}||)=cte|\pi (r_2^2-r_1^2)|\geq cte\times (r_2-r_1) \geq cte \times d, $$
 
d'ou, puisque $ v $ est reguliere,

$$ \dfrac{|a(\mu_1)-a(\mu_2)|}{|\mu_2-\mu_1|} \geq cte \dfrac{d}{|\mu_2-\mu_1|}=cte \dfrac{d(x_1,x_2)}{|v(x_1)-v(x_2)|} \geq cte >0 $$

Si on pose $ \mu_1=\mu(a) $ et $ \mu_2=\mu(b) $, on obtient le resultat.

\bigskip

\underbar { Etape 4:}

\bigskip

On pose:

$$ H(a) = \lambda \int_{\{ v>\mu(a)\} } |g'(z)|^2 e^v dx = \lambda \int_{\{ v >\mu(a) \} } e^u dx $$

$ H $ est derivable et $ H'(a)=\lambda e^{\mu(a)} $.

On a, en utilisant, l'inegalite isopermetrique, l'inegalite de Cauchy-Schwarz puis la formule de Green et la forumle de la coaire: (d'abord, on considere $ a $ tel que $ \mu(a) \in R_v $, avec $ R_v $ les valeurs regulieres de $ v $)

$$ 4 \pi a = 4 \pi a(\mu(a)) \leq \left (\int_{\{v=\mu (a)\}} |g'(z)|d\sigma_{\{v=\mu(a)\}} \right )^2  $$

$$ \leq \int_{\{v=\mu(a)\}} |\nabla v| d\sigma \int_{\{ v= \mu(a)\}} \dfrac{|g'(z)|^2}{|\nabla v|} d\sigma = $$

$$ = H(a)(-\dfrac{da}{d\mu}(\mu(a)) $$

Ce qu'on peut ecrire en considerant l'inverse de $ a $ ( car $ a $ est tel que $ \mu(a) \in R_v $):

\be -\dfrac{d\mu}{da}(a) \leq \dfrac{H(a)}{4\pi a}, \ee

La formule precedente est vraie presque partout $ a $ car l'ensemble suivant , not\'e $ C' $, $ C'= \{ a, \mu'(a) \,\,{\rm existe \,\, et \not = 0} \,\, {\rm et } \,\, \mu(a) \in C \} $ est de mesure nulle (grace a la formule de la coaire appliquee a $ \mu $), o\`u $ C $ est l'ensemble des valeurs critique de $ v $ qui est de mesure nulle par le theoreme de Sard.

1- Soit, on considere la preuve dans le Chavel de la formule de la coaire aux point r�guliers, soit on utilise la theoreme de Lebegue. Les points reguliers constituent un ouvert et la fonction $ a $ est absolument continue aux voisinages des points reguliers, elle est derivable presque partout et la ou elle derivable, elle a sa forme explicite precedent qui est non nulle. Comme elle continue est strictement decroissante, son inverse est continue est strictement decroissante, de plus  elle est derivable et sa derivee est l'inverse de la derivee de $ a $. Ce qui permet d'avoir la derniers inegalite presque partout.
 
 2- Dans ce qui suit l'ensemble $ C $ est l'ensemble des point critiques et aussi l'ensemble ou on a pas la derivee de $ a $ par le Theoreme de Lebegue (ce qu'on a dit precedement), qui est aussi de mesure nulle.
 
En effet, posons:

$$ C'= \{ a, \mu'(a) \,\,{\rm existe \,\, et \not = 0} \,\, {\rm et } \,\, \mu(a) \in C \} $$

alors,

$$ \mu(C')\subset C $$

et,

$$  |\mu(C')| \leq |C|=0.     $$ 

et,

$$ \int_{C'}(-\dfrac{d\mu}{da})(a)=\int_{ {\mathbb R} } H^0(\mu^{-1}(a) \cap C') da = \int_{ {\mathbb R} } \chi(\mu(C'))(a) da =|\mu(C')| \leq |C|=0.$$

d'ou,

$$  |C'|=0. $$

\bigskip

\underbar{In\'egalit\'e de Bol}

\bigskip

En utilisant, comme dans C. Bandle la fonction:

$$ P(a) = aH'(a)-H(a)+\dfrac{1}{8\pi}H(a)^2. $$

Alors,

$$ P'(a)=aH''(a)+\dfrac{1}{4\pi}H(a)H'(a)=\lambda a e^{\mu(a)}(\dfrac{d\mu}{da}+\dfrac{H(a)}{4\pi a})\geq 0,$$ 

et le fait que, grace a l'in\'egalit\'e isoperim\'etrique:

$$ 4\pi a(0) = 4 \pi\int_{B_1} |g'(z)|^2 dx \leq \left ( \int_{\partial B_1} e^{u/2} d\sigma \right )^2 $$

On a,

$$ \left (\int_{\partial B_1} e^{u/2} d\sigma \right)^2 \geq \left ( 4\pi - \dfrac{\lambda}{2} \int_{B_1} e^u dx \right ) \left (\int_{B_1} e^u dx \right ). $$

\bigskip

\underbar {In\'egalit\'e de Huber:}

\bigskip

Soit $ g(x,y) $ la fonction de Green du Laplacien sur $ D $, $ h $ une fonction harmonique reguliere sur $ \bar D $.

\bigskip

On pose:

$$ a(g)=\int_{\{ x, g(x,y)>g\}} e^h = \int_{ \{x, g(x,y) >g\}} |g'(z)|^2 dx, \,\, a(0)=A=\int_D e^h dx   $$

$ \partial g $ est holomorphe et est strictement positive au voisinage de la singularite $ y $ ( car $ \partial g \equiv \nabla g $ est en $ 1/r, r=|x-y| $), et sur $ \partial D $ par le principe du maximum, donc, d'apres le principe des zeros isoles, les zeros de $ \partial g $ sont isoles et comme elle ne s'annule pas au bord, ses zeros sont en nombre fini. D'ou le fait que les lignes de niveaux de $ g $ sont de mesures nulles.

\bigskip

Comme $ \nabla g $ a un nombre fini de zeros, $ a(g) $ est absolument continue et comme $ g $ est continue, $ a(g) $ est strictement decroissante, elle admet un inverse $ g(a) $ strictement decroisssante et continue.

\bigskip

Comme precedemment, $ g(a) $ est localement lipschtizienne et verifie (apres avoir utiliser, la formule de Stokes dans la couronne delimitee par la ligne de niveau $ \{g=g(a)\} $ et $ \partial D $,

$$ 4\pi a = 4\pi a(g(a)) \leq \left ( \int_{\{ g= g(a)\}} e^{h/2} d\sigma \right )^2 \leq \int_{\{g=g(a)\}} |\nabla g| d\sigma \int_{\{ g= g(a)\}} \dfrac{e^h}{|\nabla g|} d\sigma, $$

$$ \int_{\{ g=g(a)\}} |\nabla g|d\sigma = \int_{\partial D} \partial_{\nu} g d\sigma_{\partial D} =1, $$

$$ \int_{\{g=g(a)\} }\dfrac{e^h}{|\nabla g|} d\sigma = -\dfrac{da}{dg} (g(a)), $$

d'ou,

\be -g'(a) \leq \dfrac{1}{4\pi a}, \,\,\forall \, a, g(a) \in R_g, \ee

avec $ R_g $ l'ensemble des valeurs regulieres de $ g $.(c'est le complementaire de l'ensemble des valeurs critiques de $ g $ qui est fini car $ \nabla g\equiv \partial g $ a un nombre fini de zeros).

d'ou en integrant entre $ A=a(0) $ et $ a $,

\be g(a) \leq \dfrac{1}{4\pi} \log \dfrac{A}{a}.\ee

On pose,

$$ H(a)=\int_{ \{x, g(x,y) > g(a)\}} e^{h+2\alpha g(x,y)} dx, $$

Alors, $ H $ est derivable et;

$$ H'(a)=e^{2\alpha g(a)}, \,\, H(0)=0 ,\,\, H(A)=H(a(0))=\int_D e^{h+2\alpha g(x,y)} dx,$$

d'ou,

$$ H(A)=\int_D e^{h+2\alpha g(x,y)} dx=\int_0^A e^{2\alpha g(a)} da \leq \int_0^A e^{ \dfrac{ \alpha}{2\pi}  \log \frac{A}{a} } da = \dfrac{2\pi}{2\pi-\alpha} A, $$

Donc, si on utilise l'inegalite isoperimetrique de Nehari

\be 2(2\pi-\alpha)\int_D e^{h+2\alpha g(x,y)} dx \leq 4\pi A=4\pi\int_D e^h dx \leq \left (\int_{\partial D} e^{h/2} d\sigma \right )^2 \ee

On pose, en considerant $ \mu_2 $ une mesure positive:

$$ p(x)=\int_D g(x,y)2 d\mu_2, g^n(x,y)=\min \{n,g(x,y)\} $$ 

Soit, $ D_{\epsilon} =D-D^0_{\epsilon} $ avec $ D^0_{\epsilon} $ un voisinage du bord de $ D $ qui tend (en mesure) vers  $ 0 $ quand $ \epsilon  $ tend vers 0. ($\alpha_{\epsilon}=\mu_2(D_{\epsilon}) $).

Alors la fonction $ g^n $ est continue sur $ D \times D_{\epsilon} $, elle est uniformement continue. On peut ecrire cela comme suit:

$$ \forall \epsilon >0 \,\, \exists \,\, \delta >0, \,\, \forall \,\, x,y, y'\,\, d(y, y')< \delta \,\,\Rightarrow \,\,\, |g^n(x,y)-g^n(x,y')|\leq \epsilon $$

On recouvre le compact $ \bar D_{\epsilon} $ par un nombre fini d'ensembles de diametres $ < \delta $, centres en des points $ y_1, \ldots, y_{k_n} $, notes $ U_1, \ldots, U_{k_n} $.

On obtient:

$$ \int_{D_{\epsilon}} g^n(x,y) dy \leq \sum_{k=1}^{k_n}|U_k| g^n(x,y_k)+\epsilon \mu_2(D_{\epsilon}),  $$

On pose,

$$ |U_k|=\nu_k|D_{\epsilon}|=\nu_k \alpha_{\epsilon}, \,\, k=1,\ldots, k_n, \,\, {\rm avec} \,\, 0 < \nu_k < 1 \,\, {\rm et} \,\, \sum_{k=1}^{k_n} \nu_k = 1. $$

On obtient, en utilisant l'inegalite de Holder et ce qui precede:

$$ \int_D e^{h+p(x)} dx \leq e^{\epsilon \mu_2(D_{\epsilon})}\int_D e^{\sum_{k=1}^{k_n} \nu_k(h+2\alpha_{\epsilon} g(x,y_k))} dx  $$

$$ \leq e^{\epsilon \mu_2(D_{\epsilon})} \Pi_{k=1}^{k_n} \left ( \int_D e^{h+2\alpha_{\epsilon} g(x,y_k)} dx \right)^{\nu_k} \leq \dfrac{2 \pi}{2\pi-\alpha_{\epsilon}} A, $$

En faisant tendre $ \epsilon  $ vers  $ 0 $ puis $ n $ vers l'infini, en utilisant le lemme de Fatou on obtient:

\be 2(2\pi-\alpha)\int_D e^{h+\int_D g(x,y)2d\mu} \leq 4\pi A = 4\pi \int_D e^h dx \leq \left (\int_{\partial D} e^{h/2}d\sigma \right)^2, \ee

avec, $ \alpha =\mu_2(D), \,\, \mu=\mu_2-\mu_1 $

\bigskip

On applique cela avec, 

$$ \mu(B)=\int_B (K-K_0)e^u dx, \,\, \Rightarrow \,\, \mu_2(B)=\int_{\{ B, K >K_0 \} } (K-K_0) e^u dx  $$

et 

$$ \mu_1(B)=-\int_{\{B, K<K_0 \} }(K-K_0) e^u dx, $$

On pose alors:

$$ p(x)=\int_D g_D(x,y)2(K-K_0)e^u dx, \,\, q(x)=\int_B g(x,y)2(K-K_0)e^u dx, $$
 
ou $ B \subset D $.

\bigskip

On ecrit, puisque $ p-q $ est harmonique dans $ B $ et $ q=0 $ sur $ \partial B $,

$$ 2(2\pi-\mu_2(B)) \int_B e^{h+p} dx = 2(2\pi-\mu_2(B)) \int_B e^{h+(p-q)+q} dx \leq \left ( \int_{\partial B} e^{(h+p)/2} d\sigma \right)^2, $$

On ecrit alors,

\be 2(2\pi-\mu_2(D)) \int_{D_{\mu(a)}} e^{h+p} dx \leq \left ( \int_{\partial D_{\mu(a)}} e^{(h+p)/2} d\sigma \right)^2, \ee

Soit $ a $ tel que $ \mu(a) $ est valeur reguliere de $ v $, alors $ D_{\mu(a)}=\{ v>\mu(a)\} $ est une sous variete. On note $ \Omega_j $ les composantes connexes de $ D_{\mu(a)} $, alors $ \bar \Omega_i \cap \bar \Omega_j = \emptyset $ pour $ i \not = j $, pour le voir, on suppose que $ x_0 \in \bar \Omega_i \cap \bar \Omega_j $, alors en considerant un vecteur tangent $ ( v(\gamma (t))=\mu(a)) $et un vecteur rentrant $ (v(x_0 \pm h e_1) \geq \mu(a)=v(x_0)) $, on a $ \nabla v(x_0)=0 $ ce n'est pas possible car $ \mu(a) $ est valeur reguliere de $ v $.

On ecrit alors ($ q_j = 0 $ sur $ \partial \Omega_j $):

$$ \int_{D_{\mu(a)}} e^{h+p} dx = \sum_j \int_{\Omega_j} e^{h+(p-q_j)+q_j} dx \leq \sum_j \dfrac{1}{2(2\pi-\mu_2(\Omega_j))}(\int_{\partial \Omega_j} e^{(h+p-q_j)/2} d\sigma )^2 $$

$$ \leq \sum_j \dfrac{1}{2(2\pi-\mu_2(D))}(\int_{\partial \Omega_j} e^{(h+p)/2} d\sigma)^2 \leq \dfrac{1}{2(2\pi-\mu_2(D))}(\sum_j\int_{\partial \Omega_j} e^{(h+p)/2} d\sigma )^2 $$

$$ = \dfrac{1}{2(2\pi-\mu_2(D))}\left ( \int_{\partial D_{\mu(a)}} e^{(h+p)/2}\right)^2, $$

Maintenant, si on suppose que sur $ 0 \in D $, $ D $ la boule de rayon 1, qu'on a un terme $ p_0 $ tel que:

$$ p_0(x)=-2\alpha \log |x| = \int_D g(x,y) d\mu_0 =\int_D g(x,y) d \delta_0, $$

alors avec le proced\'e precedent, on a, si $ 0 \in B $:

$$ 2(2\pi-\mu_2(B)-\mu_0(B)) \int_B e^{h+p+p_0} dx = 2(2\pi-\mu_2(B)-2\pi \alpha) \int_B e^{h+(p_0-q_0)+(p-q)+q+q_0} dx  $$

$$ \leq  \left ( \int_{\partial B} e^{(h+p+p_0)/2} d\sigma \right)^2 .$$

\bigskip

\bigskip

\underbar{In\'egalit\'e d'Alexandrov}

\smallskip

Comme dans C. Bandle, on obtient, en posant,

$$ \tilde u= u-p,\,\, \tilde v= \tilde u-\tilde h,\,\, \lambda= 2K_0\geq 0 $$

(Si $ K_0=0 $, c'est l'in\'egalit\'e de Huber, on peut alors supposer que $ K_0 >0 $.)

et,

$$ \tilde a(\mu)=\int_{\{ \tilde v > \mu\}} |\tilde g'(z)|^2e^p dx, \,\,\, H(a) = \lambda \int_{\{\tilde v> \mu(a)\}} |\tilde g'(z)|^2e^p e^{\tilde v} dx, $$

et,

$$ P(a) = aH'(a)-H(a)+\dfrac{1}{4 \alpha}H(a)^2. $$

$$ \alpha = 2\pi-\mu_2(D)=2\pi-\int_{\{x\in D, K >K_0\}} (K-K_0)e^u, $$

alors,

\be H'(a)=\lambda e^{\mu(a)}, \,\, {\rm et} \,\, -\dfrac{d\mu}{da}(a) \leq \dfrac{H(a)}{2\alpha a}, \ee

et,

$$ P'(a)=aH''(a)+\dfrac{1}{2 \alpha }H(a)H'(a)=\lambda a e^{\mu(a)}(\dfrac{d\mu}{da}+\dfrac{H(a)}{2 \alpha a})\geq 0,$$ 

et le fait que, grace a l'in\'egalit\'e de Huber:

$$ \left ( \int_{\partial D} e^{u/2} d\sigma \right )^2=\left (\int_{\partial D} e^{(\tilde h+p)/2} d\sigma \right )^2 \geq 2(2\pi-\mu_2(D))\int_{D} e^{\tilde h+p} dx = 2\alpha a(0), $$

On a:

$$ \left ( \int_{\partial D}e^{u/2} d\sigma \right )^2 \geq \left (2\alpha -\dfrac{\lambda}{2} \int_D e^u dx \right ) \int_D e^u dx, $$

\bigskip

\underbar{La methode de Bol-Fiala, les geodesiques paralleles:}

\bigskip

La construction que fait Fiala, pour prouver l'inegalite isoperimetrique,se base sur une definition des coordonnees geodesiques paralleles. Il considere une courbe $ C $ delimitant un domaine simplement connexe. Cette courbe $ C $ est supposee reguliere (analytique) dans un espace de Rieman complet ( Hopf -Rinow) puis s'occupe de definir des points extremaux, points ou les geodesiques ne sont plus minimisantes, points focaux, afocaux, qui sont en nombre fini sur tout ensemble born\'e. A partir de la, on peut parcourir notre espace  partir de la courbe, en etant  une distance r de celle-ci et ne rencontrer qu'un nombre fini de points extremes stationnaires.

\smallskip

Si on note $ L(p) $ la longeur de la vraie parallele alors le but est de prouver que:

\be  \dfrac{dL(p)}{dp}  \leq  \int_0^{L(C)} k(q) dq-C(p), \,\, si \,\, p >0 \ee

\be  \dfrac{dL(p)}{dp}  \geq  \int_0^{L(C)} k(q) dq + C(p), \,\, si  \,\, p_{min} < p < 0 \ee

ou, $ C(p) $ est l'integrale de la courbure totale sur le domaine compris entre la courbe $ C $ et la vraie parallele.

\smallskip

Ceci, se fait grace  l'expression de $ L(p) $ en coordonnes geodesiques paralleles (qui sont semblables aux coordonnees de Fermi) et la formule de Gauss-Bonnet (ici, la caracteristique d'Euler-Poincare est 1, car $ C $ delimite un domaine simplement connexe $ F $. Si on note $ K $ la courbure de Gauss, on obtient:

\be 2\pi =  \int_0^{L(C)} k(q) dq+ \int_F K(x) dx \ee

Les deux inegalites precedentes, peuvent etre ecrites sous la forme suivante:

\be  \dfrac{dL(p)}{dp}  \leq  2\pi - \int_{F_p} K(x) dx, \,\, si \,\, p >0 \ee

\be  \dfrac{dL(p)}{dp}  \geq  2\pi - \int_{F_p} K(x) dx, \,\, si  \,\, p_{min} < p < 0 \ee

Avec $ F_p, p >0 $, le domaine totale, union de celui delimite par $ C $ et celui entre $ C $ et la vraie parallele. De meme pour $ p <0 $.

\smallskip

Notons que Hartman a generalise ce procede pour des courbes lipschitziennes. Les proprietes precedentes sont vraies presque partout en $ p $.

\end{document}